\documentclass[titlepage,11pt]{article}
\usepackage{latexsym}
\usepackage{amsfonts}
\usepackage{amsmath}
\usepackage{mathtools}
\usepackage{tikz}
\usepackage{comment}
\usetikzlibrary{shapes.geometric}
\newcommand\blackslug{\hbox{\hskip 1pt \vrule width 4pt height 8pt depth 1.5pt
        \hskip 1pt}}
\newcommand\bbox{\hfill \quad \blackslug \bigbreak}

\newcommand{\mab}{\mathbb}
\def\DD{\hbox{-}}
\def\CC{\hbox{-}\cdots\hbox{-}}
\def\LL{,\ldots,}
\newcommand{\vare}{\varepsilon}

\newcommand{\cupcup}{\cup \cdots\cup}

\def\polylog{\operatorname{polylog}}

\DeclarePairedDelimiter\abs{\lvert}{\rvert}%

\title{Trees and near-linear stable sets}
\author{
Tung Nguyen\thanks{Supported by AFOSR grant
FA9550-22-1-0234 and by NSF grant DMS-2154169.}\\
Princeton University,\\ Princeton, NJ 08544, USA
\and
Alex Scott\thanks{Supported by EPSRC grant EP/X013642/1}\\
University of Oxford, \\
Oxford, UK
\and
Paul Seymour\thanks{Supported by AFOSR grant
FA9550-22-1-0234 and by NSF grant DMS-2154169.}\\
Princeton University,\\ Princeton, NJ 08544, USA}

\date{May 10, 2024; revised \today}

\newtheorem{thm}{}[section]

\newcommand{\Proof}{\noindent{\bf Proof.}\ \ }

\begin{document}
\maketitle
\begin{abstract}
When $H$ is a forest, the Gy\'arf\'as-Sumner conjecture implies that every graph $G$  with no induced subgraph isomorphic to $H$ and 
with bounded clique number 
has a stable set of linear size. We cannot prove that, but we prove that every such graph $G$ has a stable set of size $|G|^{1-o(1)}$.
If $H$ is not a forest, there need not be such a stable set.

Second, we prove that when $H$ is a ``multibroom'', there {\em is} a stable set of linear size. 
As a consequence, we deduce that all multibrooms satisfy a ``fractional colouring''
version of the Gy\'arf\'as-Sumner conjecture.

Finally, we discuss extensions of our results to the multicolour setting.

\end{abstract}

\section{Introduction}

Graphs in this paper are finite, and have no loops or parallel edges. For a graph $G$, $|G|$ denotes the number of vertices of $G$,
$\chi(G)$ is its chromatic number, and $\omega(G)$ and $\alpha(G)$ denote the sizes of its largest clique and stable set 
respectively. If $G,H$ are graphs, $G$ is {\em $H$-free} if no induced subgraph of $G$ is isomorphic to $H$, and if $\mathcal{H}$
is a set of graphs, $G$ is {\em $\mathcal{H}$-free} if $G$ is $H$-free for each $H\in \mathcal{H}$. 

We are interested in how large $\alpha(G)$ must be, when $\omega(G)$ is bounded. 
For each integer $k\ge 1$, every $K_{k+1}$-free graph $G$ satisfies $\alpha(G)\ge \Omega(|G|^{1/k})$,
and there are $K_{k+1}$-free graphs $G$ such that $\alpha(G)\le |G|^{2/(k+2)}\log |G|$ (Spencer~\cite{spencer}), so at least 
we have a rough idea of the right answer.
But what happens if we fix a graph $H$ and ask the same question for $\{H,K_{k+1}\}$-free graphs $G$? 
If $H$ is a forest this question is particularly interesting, and is the focus of this paper. 

The Gy\'arf\'as-Sumner conjecture~(\cite{gyarfastree, sumner}; see also \cite{gyarfasperfect}) says:
\begin{thm}\label{GSconj}
For every forest $H$ and every integer $k\ge 1$, there exists $c\ge 1$ such that $\chi(G)\le c$ for every $\{H,K_{k+1}\}$-free graph $G$.
\end{thm}
In particular, 
if the Gy\'arf\'as-Sumner conjecture is true, then every $\{H,K_{k+1}\}$-free graph $G$
has a stable set of size at least $|G|/c$.
This is dramatically different from what happens when only $K_{k+1}$ is excluded,
and for such graphs $G$ we ought at least to be able to prove a lower bound on $\alpha(G)$ better than $|G|^{1/k}$. Until now that
had not been done,
but we will prove:
\begin{thm}\label{mainthm1}
For every forest $H$ and every integer $k\ge 1$, every $\{H,K_{k+1}\}$-free graph $G$ satisfies $\alpha(G)\ge |G|^{1-o(1)}$,
and hence has chromatic number at most $O(|G|^{o(1)})$.
\end{thm}
(The second claim follows from the first by recursive deleting maximum stable sets.)

An appealing feature of the Gy\'arf\'as-Sumner conjecture is that it is sharp, in the sense that if
$H$ is not a forest, and $k\ge 2$, then there does not exist $c$ such that 
$\chi(G)\le c$ for every $\{H,K_{k+1}\}$-free graph $G$. This is a consequence of a result of Erd\H{o}s~\cite{erdos}, who showed that 
for all $g\ge 1$ there exists 
$\vare>0$ and arbitrarily large graphs $G$ of girth at least $g$ in which every stable set has size at most $|G|^{1-\vare}$. 
Thus such graphs $G$ will be $H$-free if $H$ has a cycle of length less than $g$, and will have chromatic number more than $c$ when $|G|^{\vare}>c$. 
The same example shows that \ref{mainthm1} is sharp in the same sense: if $H$ is not a forest, and $k\ge 2$,
then not every $\{H,K_{k+1}\}$-free graph $G$ satisfies $\alpha(G)\ge |G|^{1-o(1)}$. 

To prove \ref{mainthm1}, we could assume that $H$ is a tree, and hence has a radius. (The {\em radius} of a connected graph is the 
minimum over $u\in V(G)$ of the maximum over $v\in V(G)$ of the distance between $u,v$.) Here is a more exact version of \ref{mainthm1}
(logarithms in this paper are to base two):

\begin{thm}\label{mainthm2}
Let $k,r\ge 2$ be integers, let $q:=(r-1)(k-1)$, and let $T$ be a tree of radius at most $r$. Then there exists 
$b>0$ such that every $\{T,K_{k+1}\}$-free graph $G$ satisfies
$$\alpha(G)\ge |G|^{1-b(\log |G|)^{-\frac1q}}.$$ 
\end{thm}
(This bound is better than $\Omega(|G|^{1-\vare})$, for any fixed $\vare>0$, but not as good as
$\Omega(|G|/\polylog(|G|)$.)

The same proof yields a stronger result, that implies \ref{mainthm2} by setting $d=|G|$:
\begin{thm}\label{sparsemainthm}
Let $k,r\ge 2$ be integers, let $q:=(r-1)(k-1)$, and let $T$ be a tree of radius at most $r$.
Then there exists $b>0$ such that for
every $d\ge 2$,
every $\{T,K_{k+1}\}$-free graph $G$ with
maximum degree at most $d$ satisfies
$$\alpha(G)\ge 2^{-b(\log d)^{1-\frac1q}}|G|.$$
\end{thm}

For some trees $T$, we can do better than this: we can prove that all $\{T,K_{k+1}\}$-free graphs have stable sets of linear size 
(as they should if the Gy\'arf\'as-Sumner conjecture is true), and that 
leads to the second main result of the paper. Let us explain.
A {\em broom} is a tree obtained from a path $v_0\DD v_1\CC v_{\ell}$ (where $\ell\ge 1$) by adding some number of new vertices,
each adjacent to $v_{\ell}$.  We call $v_0$ the {\em root} of the broom, and $\ell$ is its {\em length}. If $\ell_1\LL \ell_n\ge 1$,
an {\em $(\ell_1\LL \ell_n)$-multibroom} is a graph obtained from brooms of lengths $\ell_1\LL \ell_n$ by identifying their roots, and a 
{\em multibroom} is a $(\ell_1\LL \ell_n)$-multibroom for some choice of $\ell_1\LL \ell_n$.

It is not known whether all multibrooms satisfy the Gy\'arf\'as-Sumner conjecture; indeed, the only multibrooms
currently known to do so are:
\begin{itemize}
\item $(1\LL 1,2\LL 2)$-multibrooms (due to Scott and Seymour~\cite{newbrooms},
unifying earlier theorems of Gy\'arf\'as~\cite{gyarfastree}, Kierstead and Penrice~\cite{kiersteadpenrice}, and Kierstead and Zhu~\cite{kiersteadzhu});
\item $(\ell_1,\ell_2)$-multibrooms for all $\ell_1,\ell_2$, and more generally, multibrooms made by identifying the roots of several 
brooms, where all but one of the brooms is a path
(by Chudnovsky, Scott and Seymour~\cite{distantstars}; Scott \cite{scott} had previously shown the result when all brooms are paths).
\end{itemize}
Incidentally, this is almost the complete list of trees known to satisfy the Gy\'arf\'as-Sumner conjecture\footnote{The only 
other trees known to do so are those obtained by subdividing a star and adding one leaf,
attached arbitrarily~\cite{distantstars}, and those obtained from two disjoint paths by adding an edge between them~\cite{spirkl}.  
See~\cite{survey} for discussion, and for other recent work see~\cite{gsnote} and papers in the series~\cite{poly8}.}. It is known~\cite{scott} that for every tree $T$ it suffices to exclude a finite list of subdivisions of $T$.

The second main result of this paper says that all multibrooms satisfy a ``fractional'' version of \ref{GSconj}.
If $a,b\ge 1$ are integers, an {\em $(a,b)$-fractional colouring} of $G$ is a family $(A_i:1\le i\le a)$ of stable sets of 
$G$ such that 
$|\{i\in \{1\LL a\}:v\in A_i\}|\ge b$ for each vertex $v$. The {\em fractional chromatic number} 
$\chi^*(G)$ of $G$ is the minimum of $a/b$ over all pairs $a,b$ such that $G$ admits an $(a,b)$-fractional colouring. 
We will prove:
\begin{thm}\label{fracmultibrooms}
For every multibroom $T$, and every integer $k\ge 1$, there exists $c>1$ such that $\chi^*(G)\le c$ for every $\{T,K_{k+1}\}$-free graph $G$.
\end{thm}
As far as we know, this is the only case when the fractional version of the Gy\'arf\'as-Sumner conjecture has been proved and the full conjecture has not.
We will deduce this result via linear programming duality from a weighted version of the following, which strengthens \ref{mainthm2} for multibrooms to give a linear bound:
\begin{thm}\label{linmultibrooms}
For every multibroom $T$, and every integer $k\ge 1$, there exists $c>0$ such that $\alpha(G)\ge c|G|$ for every $\{T,K_{k+1}\}$-free graph $G$.
\end{thm}

In the final section, we discuss extensions of our results to the multicolour setting.

\section{Excluding a general forest}

This section contains the proof of \ref{sparsemainthm}.  We follow a strategy of iterative sparsification, that we also used in several recent papers on the Erd\H{o}s-Hajnal conjecture, for instance~\cite{density5}. Starting with a $T$-free graph $G$, we look for an induced subgraph $G'$ of $G$ such that $G'$ is significantly sparser than $G$, but not too much smaller; provided we can control the trade-off between size and density, we can then iterate the process until we end up with a large stable set.

To make this strategy work, we need to be able to carry out the sparsification step without sacrificing too much in terms of size.  We 
construct the sparse subgraph $G'$ in small steps: we repeatedly find a set $A$ of vertices that is sparse to most of the graph, and then take a union of many such sets.  The key part of the argument is to find a suitable set $A$. We do this in \ref{lem:key1} by embedding $T$ into $G$ one vertex at a time: as $G$ is $T$-free, this process must get stuck at some point, and we will use this to find the set $A$. 

Before starting the argument, we make a few comments:
\begin{itemize}
    \item The set $A$ that we obtain must satisfy a suitable sparsity condition: in fact, we get two sets $A$ and $B$ so that the 
vertices of $A$ have relatively few neighbours outside $A\cup B$.   We then repeat the argument on $G\setminus (A\cup B)$, and 
continuing in this way, we obtain a sequence of disjoint sets $A_i$ and $B_i$ such that the $A_i$ are sparse to later sets in the 
sequence.  Provided the $A_i$ are not too small and the $B_i$ are not too large, we obtain the vertices of $G'$ by taking a union of (some of) the sets $A_i$.    
    \item We may get stuck at different stages when attempting to embed $T$.  The size and sparsity of the sets $A_i$ and $B_i$ that we obtain depend on where we get stuck.  So, in fact, we end up creating sets of several different types: one of these types will occur often enough so that we can obtain $G'$ by taking all sets $A_i$ of this type.
    \item As we grow the tree, the sets $A_i$ get smaller.  By embedding the vertices in a careful order, we can ensure that our bounds depend primarily on the {\em radius} of the tree $T$, and only weakly on the number of vertices of $T$.  For this reason, and for 
convenience,
we embed our tree in a depth-first search order.
\end{itemize}

We will need the following well-known version of Ramsey's theorem (see~\cite{poly2} for example):
\begin{thm}\label{Ramsey}
For $k\ge 1$ an integer, if a graph $G$ has no stable subset of size $k$, then
$$|V(G)|\le \omega(G)^{k-1}+\omega(G)^{k-2}+\cdots+\omega(G).$$
Consequently $|V(G)|< \omega(G)^k$ if $\omega(G)>1$.
\end{thm}

Let us recall the definition of a depth-first enumeration of a tree.
Let $T$ be a tree with vertex set $\{\sigma_1\LL \sigma_{t}\}$.  The numbering $(\sigma_1\LL \sigma_{t})$ is a
{\em dfs-enumeration for $T$ rooted at $\sigma_1$} if, for each $i\in \{1\LL t-1\}$, $\sigma_{i+1}$ has a neighbour
in the path of $T$ between $\sigma_1,\sigma_i$. 
By growing a depth-first tree in $T$, it follows that every tree $T$ admits a dfs-enumeration rooted at any vertex $\rho\in V(T)$. 

Dfs-enumerations have useful properties that will be helpful in the proof of \ref{lem:sparse}.  
If
$(\sigma_1\LL \sigma_{t})$ is a dfs-enumeration for $T$, then for $1\le i\le t$ the subgraph $T_i$ of $T$
induced on $\{\sigma_1\LL \sigma_i\}$ is a tree.  
We refer to the path $P_i$ from $\sigma_1$ to $\sigma_i$ is the {\em active path} of $T_i$.  Later vertices in the dfs-enumeration 
of $T$ can only attach to vertices of $T_i$ that lie on the active path: if $j>i$ and $\sigma_j$ has a neighbour $\sigma'$ 
in $T_i$ then $\sigma'$ must lie on the active path $P_i$.  Furthermore, if $i\ge 2$ and $\sigma_i$ is a leaf then no vertex 
later in the dfs-enumeration is adjacent to $\sigma_i$.

The following two facts will be helpful:


\begin{thm}\label{lem:sparse}
        Let $G$ be a graph with clique number at most $k$, where $k\ge 1$ is an integer, 
        and let $n_0\LL n_k\ge 1$ be integers with $n_0\le |G|$.
        Then there exist $p\in \{1,\dots,k\}$ and an induced subgraph $H$ of $G$ with at least $n_{p-1}$ vertices and maximum 
degree less than $n_p$.
\end{thm}
\Proof  
We prove this by induction on $k$.  For $k=1$ the statement is immediate as then $G$ has no edges and we can take $H=G$.  So suppose $k>1$ and we have proved the statement for smaller $k$.  If $G$ has maximum degree less than $n_1$ then we can can take $H=G$.  Otherwise, pick a vertex of maximum degree and let $G'$ be the subgraph of $G$ induced by its neighbours: $G'$ has clique number at most $k-1$, and applying the inductive statement to $G'$ with parameters $n_1,\dots,n_k$ gives the required induced subgraph $H$.~\bbox

\begin{thm}
        \label{lem:local}
        For all $d,n>0$ where $n$ is an integer, and for every graph $G$ with maximum degree less than $d$,      
there is an induced subgraph of $G$ with at least $|G|/n$ vertices and with maximum degree less than $d/n$.
\end{thm}
\Proof  Let $(V_1,\ldots,V_n)$ be a partition of $V(G)$ such that $\sum_{i=1}^n\abs{E(G[V_i])}$ is minimum.  We can choose $i$ with $|V_i|\ge |G|/n$.  Suppose that $v\in V_i$.  For $j\ne i$, $v$ has at least as many neighbours in $V_j$ as in $V_i$, or else moving $v$ to $V_j$ would contradict minimality.  Thus $v$ has at most $d_g(v)/n\le d/n$ neighbours in $V_i$. This proves \ref{lem:local}.~\bbox


A {\em copy} of $H$ in $G$ is an isomorphism $\phi$ from $H$ to an induced subgraph $\phi(H)$ of $G$.
For $A,B\subset V(G)$ and $y\in(0,1)$, we say that $A$ is {\em $y$-sparse to B} if every vertex of $A$ has at most $y|B|$ neighbours in $B$.  A key lemma in the proof of \ref{sparsemainthm} shows that if $G$ does not contain a copy of $T$ then we can find sets $A$ and $B$ such that $B$ is not too large compared to $A$, and vertices of $A$ are sparse to $G\setminus B$.  We will apply this later to prove \ref{lem:key}, which will allow us to perform the sparsification step.

Here is the key lemma.

\begin{thm} \label{lem:key1}
        Let $k, r,t\ge2$ be integers, let $q:=(k-1)(r-1)$, and let 
$y_0,y_1,\ldots,y_q\in(0,1)$ be such that $y_p\le y_{p-1}/3t$ for $1\le p\le q$.
Let $T$ be a tree with $t$ vertices and radius at most $r$,
        and let $G$ be a $T$-free graph with clique number at most $k$ and maximum degree $d$,
such that $d\ge 6t/y_{q-1}$.
Then there exist $p\in\{1,\dots,q\}$ and disjoint $A,B\subseteq V(G)$ such that
\begin{itemize}
    \item $\abs A\ge y_{p-1}d/2$ and $\abs B\le 2rtd$; 
    \item $G[A]$ has maximum degree less than $y_pd$; and
    \item every vertex in $A$ has fewer than $\max(y_pd, k^t)$ 
neighbours in $G\setminus (A\cup B)$.
\end{itemize}
\end{thm}
\Proof
The proof proceeds by attempting to grow a copy of $T$ in dfs-enumeration order: since $G$ does not contain a copy of $T$ we will get stuck at some stage, and we will use this to obtain the sets $A$ and $B$.  At each stage of the argument, we identify a copy of some subtree 
$T_s$ of $T$, and some sets $A_i$ of vertices that will allow us to grow the embedding further.  Since our vertices are chosen in 
dfs order, later vertices of the tree cannot have arbitrary neighbours in $T_s$: they can only be adjacent to vertices on the active 
path $P_s$.  So for each vertex $v_i$ of $P_s$ we keep a set of $A_i$ of vertices that we can use to add a neighbour of $v_i$ to 
our embedding at a later stage (vertices at distance $r$ from the root are handled slightly differently, as they are not adjacent 
to later vertices in the dfs-enumeration).  We need the sets $A_i$ to satisfy several restrictions: the vertices of $A_i$ should be 
adjacent to $v_i$ and no other vertex of $T_s$; each pair $A_i,A_j$ of sets must satisfy some density condition (because adding a 
vertex from $A_j$ to our embedding means that its neighbours in $A_i$ can no longer be added to the embedding and have to be removed from $A_i$, and so we will only use vertices with few neighbours in $A_i$); and it will be convenient if the sets $A_i$ are appropriately sparse, as we need to produce a set $A$ with low density. The parameters depend on the distance from the root in $T$: as we get deeper, the sets $A_i$ get smaller, and we need to adjust the numbers accordingly. 

Let us make all this more precise.  
Let $\sigma_1\in V(T)$ be a vertex with distance at most $r$ to every other vertex in $T$, and let $(\sigma_1\LL \sigma_t)$
be a dfs-enumeration of $T$ rooted at $\sigma_1$. For $1\le s\le t$ let $T_s$ be the subtree of $T$ induced on 
$\{\sigma_1\LL \sigma_s\}$.

Let $s\in \{1\LL t\}$ and let $v_1\CC v_{\ell}$ be the active path in $T_s$ (between $v_1=\sigma_1$ and $v_\ell=\sigma_s$).
Thus $\ell\le r+1$. Let $\phi$ be 
a copy of $T_s$ in $G$, with vertex set $U=\phi(V(T_s))$, and for $i=1,\ldots,\ell$, let $w_i=\phi(v_i)$.
We say that $\phi$ is  {\em good} if
there are sets $(A_i)_{1\le i\le \min(\ell,r)}$ 
and integers $(p_i)_{1\le i\le  \min(\ell, r-1)}$ such that:
\begin{itemize}
\item the $A_i$ are pairwise disjoint subsets of $V(G)\setminus U$;
\item for $i=1,\dots,\ell$, every vertex in $A_i$ is adjacent to $w_i$ and has no other neighbours in 
$U$;
\item for $i=1,\dots,\min(\ell, r-1)$, we have
$p_i\in\{1,\dots,(k-1)i\}$, and
$$\abs{A_{i}}\ge (1-\frac{s}{2t})y_{p_i-1}d,$$
and $G[A_i]$ has maximum degree less than $y_{p_i}d$; 
\item if $\min\{\ell,r\}=r$ then $A_r$ is stable and $\abs{A_r}\ge t-s$; and
\item for $1\le h < i\le \min(\ell,r)$, $A_{i}$ is $\frac{1}{2t-s}$-sparse to $A_{h}$.

\end{itemize}
We refer to the collection of $A_i$ and $p_i$ as {\em references} for $\phi$.  Note that the case $i=r$ (if it occurs) is special, as we will only use $A_r$ to embed leaves of $T$ (and we will not need to define a set $A_{r+1}$).

We can now proceed with the main argument.
So suppose, for a contradiction, that there do not exist $p\in\{1,\dots,q\}$ and $A,B\subseteq V(G)$ satisfying the lemma.

There is no good copy of $T_t=T$ in $G$, since $G$ is $T$-free. On the other hand,  
there is a good copy of $T_1$. To see this, let $v$ be a vertex of degree $d$ in $G$, and let $D$ be the set of its neighbours. 
Thus $|D|=d\ge y_1d$. 
                Applying \ref{lem:sparse} to $G[D]$ with $k$ replaced by $k-1$, and $n_p=y_{p}d$ for $0\le p\le k-1$,
                there must exist $p\in\{1,\dots,k-1\}$ and $A\subseteq D$ such that 
$\abs A\ge y_{p-1}d\ge (1-\frac{1}{2t})y_{p-1}d$ and $G[A]$ has maximum degree less than $y_{p}d$. 
The map sending $\rho$ to $v$ is therefore a good copy of $T_1$.

Consequently there is a maximum value of $s\in \{1\LL t\}$ such that there is a good copy of $T_s$ in $G$, and $s$ is less than $t$.
Let $\phi$ be a good copy of $T_s$,
let $v_1\CC v_\ell$ be the active path $P_s$ of $T_s$ (where $v_1=\sigma_1$ and $v_\ell=\sigma_s$), and let$(A_i)_{1\le i\le \min(\ell,r)}$ and $(p_i)_{1\le i\le  \min(\ell, r-1)}$
be references for $\phi$.

Since we have a dfs-enumeration, $\sigma_{s+1}$ is adjacent in $T$ to some vertex $v_j$ of $P_s$, where $j\le r$.  We will show that either we can find a set $B$ so that the pair $(A_j,B)$ satisfies the lemma, or we can extend $\phi$ by taking a suitable $x\in A_i$ as the image of $\sigma_{s+1}$; this will give a contradiction, as we have assumed that neither outcome occurs.  

We start by defining $x$ and (if necessary) $B$. There are two cases, depending on the value of $j$.

If $j<r$,
we define $x\in A_j$ and $B$ as follows. Let $B$ be the set of all vertices $v\in V(G)$ satisfying either of the following conditions:
\begin{itemize}
    \item $v$ is adjacent to or equal to some vertex of $U$; or
    \item  $v$ has more than $\frac{1}{2t-s}|A_i|$ neighbours in $A_i$ for some $i\in\{1,\dots,j\}$.
\end{itemize}
There are at most $s(d+1)$ vertices that are adjacent to or equal to some 
vertex of $U$. For $i\in\{1,\dots,j\}$, the sum of vertex degrees of vertices in $A_i$ is at most $d|A_i|$ and so there are at most $(2t-s)d$ vertices with at least 
$\frac{1}{2t-s}|A_i|$ neighbours in $A_i$. Consequently 
$$|B|\le s(d+1)+j(2t-s)d\le 2r td.$$ 
Since $A_j$ is one of the references for $\phi$, we have that
$$|A_j|\ge (1-\frac{s}{2t})y_{p_j-1}d\ge y_{p_j-1}d/2$$ 
and $G[A_j]$ has maximum degree at most 
$y_{p_j}d$.  Since, by assumption, the pair $(A_i,B)$ does not satisfy the lemma,
it follows that there is a vertex $x\in A_j$ with at least $\max(y_{p_j}d, t^k)$ neighbours in $V(G)\setminus (A_j\cup B)$. This defines $x$.

The case $j=r$ is simpler, as then $\sigma_{s+1}$ is a leaf of $T$ and so has no neighbours later in the dfs-enumeration of $T$.  In this case, we choose $x\in A_j$ arbitrarily and do not need to define $B$.

In either case, let $\phi'$ be the extension of $\phi$ to $V(T_{s+1})$ defined by $\phi'(\sigma_{s+1})=x$ (and $\phi'=\phi$ otherwise). Clearly $\phi'$ is a copy of 
$T_{s+1}$ in $G$: we claim, for a contradiction, that it is good.

To show the goodness of $\phi'$, all that remains is to find references for it.
Let $v_1'\CC v_{j+1}'$ be the active path $P_{s+1}$ of $T$ between $\sigma_1$ and $\sigma_{s+1}$;
thus, $v_i' = v_i$ for $1\le i\le j$, and $v_{j+1}'= \sigma_{s+1}$.  
Let $$A_i'=A_i\setminus (N(x)\cup \{x\})$$ 
for $1\le i\le j$ (we will define $A_{j+1}'$ later if we need it).

For $i<j$, since $x$ is $\frac{1}{2t-s}$-sparse to $A_i$, it follows that 
$$|A_i'|\ge \left(1-\frac{1}{2t-s}\right)|A_i|\ge \left(1-\frac{1}{2t-s}\right)\left(1-\frac{s}{2t}\right)y_{p_i-1}d=\left(1-\frac{s+1}{2t}\right)y_{p_i-1}d.$$
For $1\le h<i\le j$, since
$A_{i}$ is $\frac{1}{2t-s}$-sparse to $A_{h}$, and $|A_{h}'|\ge (1-\frac{1}{2t-s})|A_{h}|$, it follows that $A_{i}'$ is $\frac{1}{2t-s-1}$-sparse to $A_{h}'$.

If $j=r$, then (by the conditions for references of $\phi$) $A_j$ is stable and $|A_j|\ge t-s$.
So $A_j'$ is stable and $|A_j'|\ge t-s-1$. It follows that
$A_1'\LL A_r'$ and $p_1\LL p_{r-1}$ are references for $\phi'$, and therefore $\phi'$ is good, a contradiction.

Thus $j<r$, and so (again, by the conditions for references of $\phi$) $\abs{A_{j}}\ge (1-\frac{s}{2t})y_{p_j-1}d$, and
$G[A_j]$ has maximum degree less than $y_{p_j}d$. Hence, noting that $x\in A_j$, we have
$$|A_j'|\ge |A_j|-y_{p_j}d-1\ge \left(1-\frac{s}{2t}\right)y_{p_j-1}d-y_{p_j}d-1\ge \left(1-\frac{s+1}{2t}\right)y_{p_j-1}d,$$
since $y_{p_j-1}d/(3t)\ge y_{p_j}d$ and $y_{p_j-1}d/(6t)\ge y_{q-1}d/(6t)\ge 1$.
To complete the references for $\phi'$, all that remains is to define $A_{j+1}'$. 

As $j<r$, we have defined the set $B$ above.
Let $C=N(x)\setminus (A_j\cup B)$, so $|C|\ge \max(y_{p_j}d, t^k)$ by our choice of $x$. Each vertex in $C$
is adjacent to $x$, has no neighbour in $U$, and is $\frac{1}{2t-s}$-sparse to $A_i$
and hence $\frac{1}{2t-s-1}$-sparse to $A_i'$, for $1\le i\le j$. 

If $j=r-1$, then since $|C|\ge k^t$ and $G$ has clique number at most $k$, $C$ contains a 
stable set
$A_{j+1}'$ of size $t$ by \ref{Ramsey}; then $A_1'\LL A_r'$ and $p_1\LL p_{r-1}$ are references for $\phi'$, 
again a contradiction. Finally, if $j\le r-2$, we 
apply \ref{lem:sparse} to $G[C]$, with 
$n_p= y_{p_j+p}d$ for $0\le p\le k$: we obtain $C'\subseteq C$ and $p\in \{1\LL k\}$ such that $|C'|\ge y_{p_j+p-1}d$ 
and $G[C']$ has maximum degree at most $y_{p_j+p}d$. Let $A_{j+1}':=C'$ and $p_{j+1}':=p_j+p$.  Then 
$A_1'\LL A_{j+1}'$ and $p_1\LL p_j,p_{j+1}'$ are references for $\phi'$, a contradiction. This proves \ref{lem:key1}.~\bbox

We now use \ref{lem:key1} to show that we can carry out the sparsification step in our strategy.  The idea is to apply
\ref{lem:key1} repeatedly to give a sequence of pairs $A_i, B_i$ such that the sets $A_i\cup B_i$ are pairwise disjoint, and each pair satisfies the conclusion of the lemma.  We have to be a little careful, as \ref{lem:key1} has several different outcomes; however, one of these will occur frequently enough that we can take the union of the corresponding $A_i$ to obtain the sparse induced subgraph that we are looking for.

\begin{thm} \label{lem:key}
Let $k, r,t\ge2$ be integers, let $q:=(k-1)(r-1)$, and let
$y_0,y_1,\ldots,y_q\in(0,1)$ be such that $y_i\le y_{i-1}/3t$ for $1\le i\le q$.
Let $T$ be a tree with $t$ vertices and radius at most $r$,
        and let $G$ be a $T$-free graph with clique number at most $k$ and maximum degree at most $d$.
        Then there exist $p\in\{1,\dots,q\}$ and an induced subgraph of $G$ with at least $(20qrtk^t)^{-1}y_{p-1}\abs G$ vertices 
and maximum degree less than $y_{p}d$.
\end{thm}
\Proof
Let $n\ge0$ be maximal such that for $1\le j\le n$ there are subsets $A_j,B_j\subseteq V(G)$ and an integer 
$p_j\in\{1,\dots,q\}$, 
such that $A_1,\ldots,A_n,B_1,\ldots,B_n$ are pairwise disjoint, and for $1\le j\le n$:
\begin{itemize}
\item  $\abs{A_j}>(4rt)^{-1}y_{p_j-1}\abs{B_j}$;
\item $G[A_j]$ has maximum degree less than $y_{p_j}d$; and 
\item let $C_j:=A_1\cup\cdots\cup A_j$ and $D_j:=B_1\cup\cdots\cup B_j$; then every vertex in $A_j$ has at most 
$\max(y_{p_j}d,k^t)$ neighbours in $V(G)\setminus(C_{j}\cup D_j)$.
\end{itemize}
We claim:
\\
\\
\noindent(1) {\em \label{claim:half} $C_n\cup D_n=V(G)$.}
\\
\\
                Suppose not, and let $F:=G\setminus(C_n\cup D_n)$.
Let $d'$ be the maximum degree of $F$, and suppose first that $d'\ge 6ty_{q-1}^{-1}$. 
Then by \ref{lem:key1} applied to $F$, there exist $p\in\{1,\dots,q\}$ and disjoint $A,B\subseteq V(F)$ with $\abs{A}\ge  y_{p-1}d'/2$
and $\abs B< 2rd't$ such that $G[A]$ has maximum degree less than $y_pd'$, and every vertex in $A$ has fewer than
                $\max(y_{p}d',k^t)\le\max(y_pd,k^t)$ neighbours in 
$$V(F)\setminus (A\cup B)=V(G)\setminus(C_n\cup D_n\cup A\cup B).$$
                In particular $\abs A\ge y_{p-1}d'/2\ge (4rt)^{-1}y_{p-1}\abs B$,
                and $d'\le d$, so taking $A_{n+1}:=A$ and $B_{n+1}:=B$ and $p_{n+1}:=p$ would contradict the maximality of $n$.

This implies that $d'< 6ty_{q-1}^{-1}$. But as $F$ has maximum degree $d'$, it follows that $\chi(F)\le d'+1$, and so $F$ has a stable set $A$ of size at least 
$(d'+1)^{-1}\abs F$. Let $B:=V(F)\setminus A$.  Then $|B|\le \frac{d'}{d'+1}|F|\le d'|A|$, and so 
$$|A|\ge (6t)^{-1}y_{q-1}|B|\ge (4rt)^{-1}y_{q-1}\abs{B},$$
and setting $A_{n+1}:=A$ and $B_{n+1}:=B$ and $p_{n+1} = q$ would contradict the maximality of $n$.
This proves (1).
        
\bigskip

Thus we have $C_n\cup D_n=V(G)$.
For every $p\in\{1,\dots,q\}$, let $J_p:=\{j\in \{1\LL n\}:p_j=p\}$;
then 
$$V(G)=\bigcup_{p=1}^q\bigcup_{j\in J_p}A_j\cup B_j,$$
and so there exists $p\in\{1,\dots,q\}$ with 
$$\bigcup_{j\in J_p}\abs{A_j\cup B_j}\ge q^{-1}|G|.$$
        Let $C:=\bigcup_{j\in J_p}A_j$ and $D:=\bigcup_{j\in J_p}B_j$;
        then $\abs C\ge (4rt)^{-1}y_{p-1}\abs D$ since $\abs{A_j}\ge (4rt)^{-1}y_{p-1}\abs {B_j}$ for all $j\in J_p$.
       Hence $|D|\le 4rty_{p-1}^{-1}|C|$, and so $|C\cup D|\le (1+4rty_{p-1}^{-1})|C|\le 5rty_{p-1}^{-1}|C|$
(because $y_{p-1}\le 1$ and $r,t\ge 1$).
We deduce that 
        $$\abs C\ge (5rt)^{-1}y_{p-1}\abs{C\cup D}\ge (5qrt)^{-1}y_{p-1}\abs G.$$
        Now, for each $j\in J_i$, at most 
$$y_pd|A_j|/2 + \max(y_pd,k^t)|A_j|\le \frac32\max(y_pd,k^t)|A_j|$$ 
edges have both ends in $\bigcup_{j'\in J_i,j'\ge j}A_{j'}$ and at least one end in $A_j$;
        and so $G[C]$ has at most $\frac32\max(y_{p}d,k^t)\abs C$ edges.

Tur\'an's theorem~\cite{turan} implies that for all $n,c$,  every $n$-vertex graph with average degree at most $c$ has a stable set of size at least $n/(c+1)$.  Thus if $G[C]$ has at most $\frac32k^t\abs C$ edges, then by Tur\'an's theorem, $G[C]$ has a stable set $S$ with 
$$\abs S\ge \left(3k^t+1\right)^{-1}\abs C\ge \left(4k^t\right)^{-1}\abs C\ge \left(20qrtk^t\right)^{-1}y_{p-1}\abs G$$
and the theorem holds. 
                Otherwise, $G[C]$ has at most $\frac32y_pd |C|$ edges. In this case
                let $S'\subseteq C$ be the set of vertices of degree at most $6y_{p}d$ in $G[C]$.
                Then $\abs{C\setminus S'}\le \frac12\abs C$ and so $|S'|\ge \frac12\abs C$.
Since $G[S']$ has maximum degree at most  $6y_{p}d$, there exists $S\subseteq S'$ with $|S|\ge |S'|/6$ 
such that $G[S]$ has maximum degree at most $y_{p}d$, by \ref{lem:local}.
So
$$|S|\ge |S'|/6\ge \abs C/12\ge (60qrt)^{-1}y_{p-1}\abs G\ge \left(20qrtk^t\right)^{-1}y_{p-1}\abs G,$$
and again the theorem holds.
        This proves \ref{lem:key}.~\bbox


Finally, by iterating the sparsification given by \ref{lem:key}, we obtain \ref{sparsemainthm}, which we restate:
\begin{thm}\label{sparsemainthmagain}
Let $k,r\ge 2$ be integers, let $q:=(r-1)(k-1)$, and let $T$ be a tree of radius at most $r$. 
Then there exists $b>0$ such that for
every $d\ge 2$, 
every $\{T,K_{k+1}\}$-free graph $G$ with
maximum degree at most $d$ satisfies
$$\alpha(G)\ge 2^{-b(\log d)^{1-\frac1q}}|G|.$$
\end{thm}
\Proof
Let $t:=\abs T$ and $c=20qrtk^t$.  We will prove that the theorem holds with $b:=\log(4c^2)$. 

Let $d\ge 2$, and let $G$ be a $\{T,K_{k+1}\}$-free graph with maximum degree at most $d$.
Let $x:=(\log d)^{1/q}$.
        Suppose first that $x\le b/2$, and so $b(\log d)^{-1/q}=b/x\ge 2$. 
Since $G$ has maximum degree at most $d$, it has a stable set of size at least 
$$|G|/(d+1) \ge d^{-2}|G|\ge d^{-b\left(\log d\right)^{-\frac1q}}|G|$$
(since $b,d\ge 2$),
and so the theorem holds. Hence we may assume that $x\ge b/2$.

        Let $G_0:=G$ and let $d_0$ be the maximum degree of $G_0$; then $d_0\le d$.
        For $i=0,\dots,q$, let $y_i:=2^{-x^i}$;
then $y_{i-1}/y_{i}=y_{i-1}^{1-x}\ge 2^{x-1}\ge 2^{b/2-1}\ge  3t$ for all $i\in\{1,\dots,q\}$.
Inductively, \ref{lem:key} implies that for each $j\ge 1$, there exists $p_j\in \{1,\dots,q\}$ and an induced subgraph $G_j$ 
of $G_{j-1}$ with $\abs{G_j}\ge c^{-1}y_{p_j-1}\abs{G_{j-1}}$ 
and maximum degree $d_j\le y_{p_j}d_{j-1}$. 

Choose $n\ge 0$ minimal such that $d_n<1$.
        The minimality of $n$ implies
        $$2^{-x^q}=1/d \le d_{n-1}/d\le y_{p_1}y_{p_2}\cdots y_{p_{n-1}}\le y_1^{n-1}\le 2^{-x(n-1)}.$$
        Consequently $x(n-1)\le x^q$, and so $n\le x^{q-1}+1\le 2x^{q-1}=2(\log d)/x$.
        Let 
$$y:=y_{p_1-1}y_{p_2-1}\cdots y_{p_n-1};$$ 
then $y^x=y_{p_1}y_{p_2}\cdots y_{p_n}$
        and 
$$1\le \left(y_{p_1}y_{p_2}\cdots y_{p_{n-1}}\right)d=\left(y^x/y_{p_n}\right)d\le y^xd^2;$$ 
and so $y\ge d^{-2/x}$.
        Thus
                $$\abs{G_n}\ge c^{-n}y\abs{G} \ge c^{-2(\log d)/x }d^{-2/x}\abs G
                 = d^{-2(\log c)/x-2/x}\abs G
                 =d^{-b/x} \abs G=d^{-b(\log d)^{-1/q}}\abs G.$$
        Since $G_n$ has maximum degree less than $1$, this proves \ref{sparsemainthmagain}.~\bbox

\section{Multibrooms and linear-sized stable sets}

Now we move on to our second main theorem, \ref{fracmultibrooms}, which we will prove in this section.
For two integers $\ell,m\ge 0$, let us say
an {\em $(\ell,m)$-broom} is a tree obtained from a path of length $\ell$, with ends $a,b$ say,
by adding $m$ new vertices each adjacent to $b$.  We call $a$ the {\em root} of the broom.
If $G$ is a graph and $S\subseteq V(G)$, we denote by $N_G(S)$ the set of vertices in $V(G)\setminus S$ with a neighbour in $S$.
A {\em weighting} on a graph $G$ is a function $w\colon V(G)\to\mab R^+$;
and for $S \subseteq V(G)$ we define $w(S):=\sum_{v\in S}w(v)$. We write $w(H)$ for $w(V(H))$ when $H$ is an induced subgraph of $G$, and define the {\em maximum weighted degree} $\Delta(G,w):=\max_{v\in V(G)}w(N_G(v))$,

The linear programming duality theorem implies that, if $H$ is a hypergraph, and each vertex of $H$ belongs to at least one hyperedge,
then the following are equal:
\begin{itemize}
\item the minimum of $\sum_{A\in H}q(A)$, over all maps $q:H\to \mathbb{R}_+$ such that, for each $v\in V(H)$,
$$\sum(q(A):A\in H \text{ and }  A\ni v)\ge 1;$$
\item the maximum of $\sum_{v\in V(H)}w(v)$, over all maps $w:V(H)\to \mathbb{R}_+$ such that, for each $A\in H$,
$$\sum(w(v):v\in V(H) \text{ and } v\in A)\le 1.$$
\end{itemize}
By applying this to the hypergraph with edge-set the set of stable sets of $G$, we deduce that \ref{fracmultibrooms} is equivalent to the following:
\begin{thm}\label{lpmulti}
Let $T$ be a multibroom and $k\ge 1$ an integer. Then there exists $c>0$ such that if $G$ is $\{T,K_{k+1}\}$-free,  
and $w$ is a weighting on $G$, then
there is a stable set $S$ of $G$ with $w(S)\ge cw(G)$.
\end{thm}

One way to try to find a stable set containing a linear fraction of the total weight, is via degeneracy. A graph is {\em $d$-degenerate} if every non-null subgraph has a vertex of degree at most $d$. (It is elementary that every $d$-degenerate graph is $(d+1)$-colourable.)
If we can find an induced subgraph $G'$ of $G$ that is $d$-degenerate (where $d$ is some constant, depending on $T$ and $k$ only)
and with total weight at least  a constant fraction of the total weight of $G$, then since $G'$ is 
$(d+1)$-colourable, it has a stable set containing at least a $1/(d+1)$ fraction of the total weight of $G'$, and so we win.
We will try to assemble such a large $d$-degenerate subgraph iteratively. Let us say a subset $X\subseteq V(G)$
is {\em $d$-degenerate in $G$} if it can be ordered $X=\{x_1\LL x_n\}$ such that for $1\le i\le n$, $x_i$ has at most $d$
neighbours in $\{x_{i+1}\LL x_n\}\cup (V(G)\setminus X)$. The advantage of this more restrictive definition is that if $X$ is $d$-degenerate in $G$,
and $X'$ is $d$-degenerate in $G\setminus X$, then $X\cup X'$ is $d$-degenerate in $G$, as is easily seen. This will allow us 
to grow a large $d$-degenerate subgraph by piecing together smaller subsets recursively.

        
We start with a lemma that we can use to grow a single branch of a multibroom.  We show that, for an arbitrary vertex $v$ in a weighted graph $G$, either we can grow a broom of any specified length rooted at $v$, or we can find two disjoint sets $X,Y$ of vertices such that $w(X)$ is at least a constant fraction of $w(Y)$ 
and $X$ has bounded degeneracy in $G\setminus(Y\cup \{v\})$. We will argue later by induction on the clique number $k$, so we assume that \ref{lpmulti} holds for induced subgraphs with smaller clique number.
\begin{thm}
\label{lem:broomexp}
Let $k,\ell\ge 1$ and $m\ge 0$ be integers, and let $d\ge 1$. 
Let $G$ be a graph with clique number at most $k$, let $w$ be a weighting of $G$, and let $v$ be a vertex of $G$. 
Suppose that, 
for every induced subgraph $G'$ of $G$ with clique number less than $k$, there is a stable set $S$ of $G'$
with $w(S)\ge w(G')/d$.  Then one of the following holds:
\begin{itemize}
    \item there is an induced $(\ell,m)$-broom in $G$ with root $v$; or
    \item there exist disjoint $X,Y\subseteq V(G)\setminus \{v\}$ with $w(X)\ge w(X\cup Y)/(d2^{2\ell})$ and $w(X\cup Y)\ge w(N_G(v))$ such that
$X$ is $k^m$-degenerate in $G\setminus (Y\cup \{v\})$.
\end{itemize}
\end{thm}
\Proof We proceed by induction on $\ell$.
Let $G,v,w$ be as in the theorem, and let $\Delta=\Delta(G,w)$ be the maximum weighted degree. 
We may delete all vertices $u\ne v$ with $w(u)=0$: at the end of an argument, either we get the desired broom, or we obtain sets $X$ and $Y$ as in the lemma, and we can add all the deleted vertices to $Y$.  Thus we may
assume that $w(u)>0$ for all $u\ne v$.

One approach to finding a broom would be to start from $v$ and explore into the graph by distance from $v$: for $i\ge0$, let $U_i$ 
be the set of vertices at distance $i$ from $v$.  If any vertex $u$ in $U_\ell$ has $k^m$ neighbours outside 
$U:=\bigcup_{i=0}^\ell U_i$, then we can construct a broom by taking a shortest path from $v$ to $u$ and finding a stable set of 
size $m$ among the neighbours of $u$ outside $U$, by \ref{Ramsey}.  Otherwise, vertices in $U$ have few neighbours in the rest of 
the graph.  However this approach is too simple: the sets $U_i$ might be very large, or might quickly exhaust the graph.  We 
therefore explore the graph more carefully.  Starting with $R_0=\{v\}$, we grow a sequence of sets $R_i$, where each set $R_i$ is a 
subset of the neighbours $L_i$ of $R_{i-1}$ in the unexplored part of the graph.  We choose $R_{i-1}$ carefully, so that $L_i$ is not too heavy.

Let $L_0=\{v\}$ and $J_0=\emptyset$. For $i=1,\dots,\ell$, we define $J_{i-1}, R_{i-1}, L_i$ inductively as follows:
\begin{itemize}
    \item Let $R_{i-1}$ be a maximal subset of $L_{i-1}$ such that $N_G(R_{i-1})\setminus (J_{i-1}\cup \{v\})$ has weight at most $2\Delta$. 
    \item $L_i=N_G(R_{i-1})\setminus (J_{i-1}\cup \{v\})$. 
    \item $J_{i}=L_1\cupcup L_{i}$.
\end{itemize}
Thus $L_1=N_G(v)$ and $w(L_1)=w(N_G(v))\le \Delta$.  Note that for any vertex $w$ of $L_i$, we can find an induced path $P_w=u_1\CC u_i$ such that $u_1=v$, $u_i=w$ and $u_j\in R_j$ for each $j$; furthermore, for $h<i$, $u_h$ has no neighbours outside $J_{i}$.
\\
\\
(1) {\em We may assume that $w(L_{j+1})< \min(w(L_{j}), \Delta)$ for some $j\in \{1\LL \ell-1\}$.}
\\
\\
Suppose first that $w(L_i)\ge \min(w(L_{i-1}),\Delta)$ for all $i\in \{2\LL \ell\}$. Since $w(N_G(v))=w(L_1)\le \Delta$, it follows that 
$w(L_{\ell})\ge w(N_G(v))$.

We claim that $w(J_{\ell-1})\le 2(\ell-1) w(L_{\ell})$. If  $w(L_{\ell})\ge \Delta$,
this is clear, 
since $w(L_i)\le 2\Delta$ by definition. If $w(L_{\ell})< \Delta$, then by applying the condition $w(L_i)\ge \min(w(L_{i-1}),\Delta)$ repeatedly, we
see that $\Delta>w(L_{\ell})\ge\dots\ge w(L_1)$.  So $w(J_{\ell-1})\le (\ell-1) w(L_{\ell})$.  In either case,
$w(J_{\ell-1})\le 2(\ell-1) w(L_{\ell})$.

Let $R_{\ell-1} = \{v_1\LL v_n\}$, numbered in arbitrary order, and for $1\le i\le n$ let $C_i$
be the set of vertices in $L_{\ell}$ that are adjacent to $v_i$ and nonadjacent to each of $v_1\LL v_{i-1}$.
Then $C_1\LL C_n$ are pairwise disjoint and have union $L_{\ell}$.  For each $i$, $G[C_i]$ has clique number 
less than $k$, as $C_i$ is contained in the neighbourhood of $v_i$; and so, by a hypothesis of the theorem, there is a stable set $S_i$ of $G[C_i]$
with $w(S_i)\ge  w(C_i)/d$. Let $S=S_1\cupcup S_n$; thus $w(S)\ge w(L_{\ell})/d$, and we have already shown that $w(L_{\ell})\ge w(N_G(v))$.


Let $Y=J_{\ell-1}\cup (L_{\ell}\setminus S)$ and,
for $1\le i\le n$, let
$$Q_i=S_{i+1}\cupcup S_n\cup (V(G)\setminus (J_{\ell-1}\cup \{v\})).$$
No vertex $u\in C_i$ has $m$ pairwise nonadjacent 
neighbours in $Q_i$, or we could find an induced $(\ell,m)$-broom by taking the path $P_w$ and adding these neighbours.
Thus $u$ has fewer than $k^m$ neighbours in $Q_i$.
Consequently $S$ is $k^m$-degenerate in 
$G\setminus (Y\cup \{v\})$. Since 
$w(S\cup Y)\le 2\ell w(L_{\ell})\le 2d\ell w(S)\le d2^{2\ell}w(S)$,
the theorem holds in this case.
This proves (1).

\bigskip
We may therefore assume that $w(L_{j+1})< \min(w(L_{j}), \Delta)$ for some $j\in \{1\LL \ell-1\}$; choose $j$ 
minimal with this property.  From the minimality of $j$, it follows that
$w(J_{j-1})\le 2(j-1)w(L_j)$, and $w(L_j)\ge w(N_G(v))$, as in the proof of (1).

Suppose first that $j=1$. Since $G[L_1]$ has clique number less than $k$, and $L_1\ne \emptyset$, 
there is a nonempty stable subset $S\subseteq L_1$ with $dw(S)\ge w(L_1)$.
Let $Y=(L_1\setminus C)\cup L_2$. Thus 
$$w(S\cup Y)\le w(L_1)+ w(L_2) \le 2w(L_1)\le 2dw(S)\le 2^{2\ell}dw(S),$$
and $S$ is $0$-degenerate and hence $k^m$-degenerate in $G\setminus (Y\cup \{v\})$. Hence the theorem holds in this case.

So we may assume that $j\ge 2$.
Choose disjoint $X,Y'\subseteq L_j$ with $X\cup Y'$ maximal such that 
$X$ is $k^m$-degenerate in $G[L_j\setminus Y']$ and $w(X\cup Y')\le 2^{2(\ell-j+1)}d w(X)$. 

Suppose that $X\cup Y'\ne L_j$. Choose $q\in L_j\setminus (X\cup Y')$, and let $p\in R_{j-1}$ be adjacent to $q$.
Let $G'=G[(L_j\setminus (X\cup Y'))\cup \{p\}]$. 
There is a path $P$ of length $j-1$ between $v,p$ with interior in $R_1\cupcup R_{j-2}$, and so 
there is no induced $(\ell-j+1,m)$-broom in $G'$ with root $p$, since its union with $P$ would be an 
induced $(\ell,m)$-broom in $G$ with root $v$. From the inductive hypothesis, applied to $G'$ and $p$, 
there are disjoint $X'',Y''\subseteq V(G')\setminus \{p\}$ with $X''\ne \emptyset$ (since $w(q)>0$) such that
$w(X''\cup Y'')\le 2^{2(\ell-j+1)}d w(X'')$ and
$X''$ is $k^m$-degenerate in $G'\setminus (Y''\cup \{p\})$. But then $X\cup X''$ is $k^m$-degenerate in 
$G[L_j\setminus(Y'\cup Y'')]$, and 
$$w(X\cup X''\cup Y'\cup Y'')\le 2^{2(\ell-j+1)}d w(X\cup X''),$$ 
contrary to the maximality of $X\cup Y'$.

Thus $X\cup Y'=L_j$. Then 
$X$ is $k^m$-degenerate in $G\setminus (Y\cup \{v\})$, where $Y=Y'\cup J_{j-1}\cup L_{j+1}$. But
\begin{align*}
w(X\cup Y')&= w(L_j)\\
w(J_{j-1})&\le 2(j-1) w(L_j)\\
w(L_{j+1})&\le w(L_j)
\end{align*}
and it follows that 
$$w(X\cup Y)\le 2jw(L_j)\le (2j)2^{2(\ell-j+1)}d w(X)\le 2^{2\ell}dw(X).$$
Since 
$$w(X)\ge (2^{2(\ell-j+1)}d)^{-1}w(X\cup Y')\ge (2^{2(\ell-j+1)}d)^{-1}w(N_G(v))\ge 2^{-2\ell}d^{-1}w(N_G(v)),$$
the theorem holds in this case.
This proves \ref{lem:broomexp}.~\bbox


Now we will apply \ref{lem:broomexp} to prove \ref{lpmulti}, which we restate:
\begin{thm}\label{lpmultiagain}
Let $T$ be a multibroom and $k\ge 1$ an integer. Then there exists $c>0$ such that if $G$ is $\{T,K_{k+1}\}$-free,
and $w$ is a weighting on $G$, then
there is a stable set $S$ of $G$ with $w(S)\ge cw(G)$.
\end{thm}
\Proof
We may assume that $k\ge 2$, and assume inductively that the theorem holds for smaller values of $k$. Thus there exists
$d\ge 1$ such that for every $\{T,K_{k}\}$-free graph $G'$, if 
$w$ is a weighting on $G'$, then
there is a stable set $S$ of $G'$ with $dw(S)\ge w(G')$.

Let $B_1,B_2,\ldots,B_b$ be vertex-disjoint brooms such that $T$ is obtained by identifying their roots.
        For each $i\in\{1\LL b\}$, let $B_i$ be an $(\ell_i,m_i)$-broom, and let $\ell=\max(\ell_1\LL \ell_b)$, and $m=\max(m_1\LL m_b)$.
        We claim that $c=(d^2 2^{2\ell+3}|T|(k^m+1))^{-1}$ satisfies the lemma.
        To see this, let $G$ be $\{T,K_{k+1}\}$-free, and let $w$ be a weighting on $G$. Choose disjoint $X,Y\subseteq V(G)$
with $X\cup Y$ maximal such that:
\begin{itemize}
\item $d^2 2^{2\ell+3}|T|w(X)\ge w(X\cup Y)$; and
\item $X$ is $k^m$-degenerate in $G\setminus Y$.
\end{itemize}
(This is possible since we may take $X=Y=\emptyset$.)

Suppose first that $X\cup Y=V(G)$. Then $d^2 2^{2\ell+3}|T|w(X)\ge w(G)$, and since $G[X]$ is $k^m$-degenerate and hence $(k^m+1)$-colourable,
it has a stable set $S$ with $w(S)\ge w(X)/(k^m+1)\ge cw(G)$, and the theorem holds. Thus, we may assume (for a contradiction)
that $X\cup Y\ne V(G)$. Let $G'=G\setminus (X\cup Y)$, let $\Delta= \Delta(G',w')$, where $w'$
is the restriction of $w$ to $V(G')$, and choose $v\in V(G')$ with 
$w(N_{G'}(v))=\Delta$.
Since $G[N_{G'}(v)]$ has clique number less than $k$, there is a stable subset $S\subseteq N_{G'}(v)$ with 
$dw(S)\ge w(N_{G'}(v))=\Delta$.
Let $Z$ be the set of vertices $z$ of $G'$ not in $S\cup \{v\}$, such that 
$w(N_{G}(z)\cap S)\ge w(S)/(2|T|)$. Let $P$ be the sum of $w(u)w(z)$ over all adjacent $u,z$ 
with $u\in S$ and $z\in Z$; then 
$$w(Z)w(S)/(2|T|)=\sum_{z\in Z}w(z)w(S)/(2|T|) \le P\le \sum_{u\in S}w(u)\Delta= w(S)\Delta,$$ 
and so $w(Z)\le 2|T|\Delta$. 
Let $Z'=N_{G'}(v)\setminus S$.
\\
\\
(1) {\em 
Let $H$ be a connected induced subgraph of $G'\setminus (Z\cup Z')$, such that $v\in V(H)$ and $|H|\le  |T|$. 
Then for all $\ell', m'$ with $1\le \ell'\le \ell$ and $0\le m'\le m$, there is an $(\ell', m')$-broom $B$ 
in $G'\setminus (Z\cup Z')$, with root $v$, with $V(B\cap H)=\{v\}$, and such that $B\cup H$ is an induced 
subgraph of $G'$.}
\\
\\
Let $W$ be the set of vertices of $G'\setminus (Z\cup Z'\cup \{v\})$ that are equal or adjacent to some vertex of $H\setminus \{v\}$. 
Thus $w(W\cup \{v\})\le |T|\Delta$. 
Since 
$w(N_{G'}(u))\le w(S)/(2|T|)$ for each vertex $u$ of
$H\setminus v$, it follows that $w(W\cap S)\le w(S)/2$. Suppose that the desired broom $B$ does not exist. 
Thus there is no $(\ell', m')$-broom with root $v$ in $G''$, where $G''=G'\setminus (Z\cup Z'\cup W)$. 
By 
\ref{lem:broomexp} applied to $G''$, and since $\ell'\le \ell$ and $m'\le m$,
there exist disjoint $X',Y''\subseteq V(G'')\setminus \{v\}$ with 
$d2^{2\ell}w(X')\ge w(X'\cup Y'')\ge w(N_{G''}(v))$ such that
$X'$ is $k^m$-degenerate in $G''\setminus (Y''\cup \{v\})=G'\setminus Y'$, where 
$$Y'=Y''\cup Z\cup Z'\cup W\cup \{v\}.$$
Consequently $X\cup X'$ is $k^m$-degenerate in $G\setminus (Y\cup Y')$.
But
$$w(N_{G''}(v))\ge w(S)-w(W\cap S)\ge w(S)/2\ge \Delta/(2d),$$
and so $d2^{2\ell}w(X')\ge \Delta/(2d)$.
Furthermore,
\begin{align*}
w(Z)&\le 2|T|\Delta\\
w(Z')&\le \Delta\\
w(W\cup \{v\})&\le |T|\Delta
\end{align*}
and so 
$$w(Z\cup Z'\cup W\cup \{v\})\le (3|T|+1)\Delta\le (3|T|+1)d^2 2^{2\ell+1}w(X').$$
Consequently 
$$w(X'\cup Y')\le d2^{2\ell}w(X')+(3|T|+1)d^2 2^{2\ell+1}w(X')\le d2^{2\ell}w(X')(1+2(3|T|+1)d )\le d^2 2^{2\ell+3}|T|w(X').$$
But $X\cup X'$ is $k^m$-degenerate in $G\setminus (Y\cup Y')$, contradicting the maximality of $X\cup Y$. This proves (1).

\bigskip
By $b$ applications of (1), it follows that $G'$ contains a copy of $T$, a contradiction. This proves \ref{lpmultiagain}.~\bbox

\section{A multicolour version}

One can look for ``multicolour'' versions of these theorems. One simple multicolour extension of the 
Gy\'arf\'as-Sumner conjecture is:
\begin{thm}\label{GSmulticolour}
{\bf Conjecture: }For every forest $H$, and all integers $k,t\ge 1$, there exists $c>0$ such that if $G$ is a graph with 
clique number at most $k$, and
the edges of $G$ can be coloured with $t$ colours in such a way that
 no monochromatic induced subgraph of $G$ is isomorphic to $H$,
then $\chi(G)\le c$.  
\end{thm}
({\em Monochromatic} means with all edges the same colour.) Calling this a ``conjecture'' is a stretch, because we have very little faith in it.
So far we cannot even prove it when $H$ is the four-vertex path (although we can if in addition $t\le 3$).

On the other hand, surprisingly to us, \ref{GSmulticolour} is  ``nearly'' true.
Our first main theorem \ref{mainthm1} can be extended to the following:

\begin{thm}\label{multicolour}
For every forest $H$, and all integers $k,t\ge 1$, if 
$G$ is a graph with 
clique number at most $k$, and
the edges of $G$ can be coloured with $t$ colours in such a way that
no monochromatic induced subgraph of $G$ is isomorphic to $H$,
then $\alpha(G)\ge |G|^{1-o(1)}$.
\end{thm}
It is straightforward to adapt the proof of \ref{mainthm1} to give this more general result, and we omit the details.
Pleasingly, this is tight, in the sense that if $\mathcal{H}$ is a finite set of $t$-edge-coloured graphs (that is, graphs $H$ with $E(H)$
partitioned into $t$ subsets labelled $1\LL t$), and we want to know whether every $t$-edge-coloured graph with bounded clique number 
that contains no member of $\mathcal{H}$ as an induced subgraph (in the natural sense) has a nearly-linear stable set, then \ref{multicolour} tells us the answer:
if and only if for $1\le i\le t$, there is a forest in $\mathcal{H}$ with all edges coloured $i$.

Another multicolour extension of \ref{GSconj} generalizes a theorem of Chudnovsky and Seymour~\cite{ChudSey}, who proved the following when 
$t=2$:
\begin{thm}\label{ChudSey}
Let $t\ge 1$ be an integer, and let $\mathcal{H}$ be a finite set of $t$-edge-coloured complete graphs. 
Assuming the truth of \ref{GSconj}, the following 
are equivalent:
\begin{itemize}
\item 
there exists $c>0$ such that for every $t$-edge-coloured
complete graph $G$ with no induced subgraph in $\mathcal{H}$, $V(G)$ is the union of $c$ monochromatic cliques
\item for every two colours $i,j\in \{1\LL t\}$, $\mathcal{H}$ contains a forest with all edges in colour $i$ and all nonedges 
in colour $j$, and $\mathcal{H}$ contains a complete multipartite graph with all edges in colour $i$ and all nonedges in colour $j$.
\end{itemize}
\end{thm}
This can be deduced from a result of~\cite{pairs} (we omit the details), although deriving the second bullet from the first 
assumes the Gy\'arf\'as-Sumner conjecture. But without assuming that, we can prove something similar:

\begin{thm}\label{multimulti}
Let $t\ge 1$ be an integer, and let $\mathcal{H}$ be a finite set of $t$-edge-coloured complete graphs. 
The following 
are equivalent:
\begin{itemize}
\item 
every $t$-edge-coloured complete graph with no induced subgraph in $\mathcal{H}$ has a monochromatic clique of size $|G|^{1-o(1)}$;
\item for every two colours $i,j\in \{1\LL t\}$, $\mathcal{H}$ contains a forest with all edges in colour $i$ and all nonedges 
in colour $j$, and $\mathcal{H}$ contains a complete multipartite graph with all edges in colour $i$ and all nonedges in colour $j$.
\end{itemize}
\end{thm}
This can easily be deduced from \ref{mainthm1} and the theorem of \cite{pairs}, and we omit the details.



\begin{thebibliography}{99}

\bibitem{distantstars} M. Chudnovsky, A. Scott and P. Seymour,
``Induced subgraphs of graphs with large chromatic number. XII.
Distant stars'', {\em J. Graph Theory} {\bf 92} (2019), 237--254,
{\tt arXiv:1711.08612}.

\bibitem{ChudSey} M. Chudnovsky and P. Seymour, ``Extending the Gy\'arf\'as-Sumner conjecture'',
{\em J. Combinatorial Theory, Ser B},
{\bf 105} (2014), 11--16.

\bibitem{pairs} M. Chudnovsky, A. Scott and P. Seymour, ``Excluding pairs of graphs'',
{\em J. Combinatorial Theory, Ser B}, {\bf 106} (2014), 15--29.


\bibitem{erdos} P. Erd\H{o}s,
``Graph theory and probability'',
{\em Canad. J. Math.} {\bf 11} (1959), 34--38.

\bibitem{gyarfastree}
A. Gy\'arf\'as, ``On Ramsey covering-numbers'',
{\em Coll. Math. Soc. J\'anos Bolyai}, in {\em Infinite and Finite Sets},
North Holland/American Elsevier, New York (1975), 10.

\bibitem{gyarfasperfect}
A. Gy\'arf\'as, ``Problems from the world surrounding perfect graphs'', 
Proceedings of the International Conference on Combinatorial Analysis and its Applications, (Pokrzywna, 1985), 
{\em Zastos. Mat.} {\bf 19} (1987), 413--441.

\bibitem{kiersteadpenrice} 
H. A. Kierstead and S. G. Penrice,
``Radius two trees specify $\chi$-bounded classes'',
{\em J. Graph Theory} {\bf 18} (1994), 119--129.

\bibitem{kiersteadzhu} H. A. Kierstead and  Y. Zhu,
``Radius three trees in graphs with large chromatic number'',
{\em SIAM J. Disc. Math.} {\bf 17} (2004), 571--581.

\bibitem{density5} T. Nguyen, A. Scott and P. Seymour, 
``Induced subgraph density. V. All paths approach
Erd\H{o}s-Hajnal'',
isubmitted for publication, {\tt arXiv:2307.15032}.


\bibitem{gsnote} T. Nguyen, A. Scott and P. Seymour,
``A note on the Gy\'arf\'as-Sumner conjecture'', 
{\em Graphs and Combinatorics} {\bf 40} (2024), article 33.


\bibitem{poly8} T. Nguyen, A. Scott and P. Seymour,
``Polynomial bounds for chromatic number. VIII. Excluding a path and a complete multipartite graph'', 
{\em J. Graph Theory} {\bf 107} (2024), 509--521, {\tt arXiv:2303.11766}.

\bibitem{scott} A. Scott,
``Induced trees in graphs of large chromatic number'', 
{\em J. Graph Theory} {\bf 24} (1997), 297--311.

\bibitem{newbrooms} A. Scott and P. Seymour,
``Induced subgraphs of graphs with large chromatic number. XIII.
New brooms'',
{\em European J. Combinatorics} {\bf 84} (2020), 103024, {\tt arXiv:1807.03768}.

\bibitem{survey}
A. Scott and P. Seymour, 
``A survey of $\chi$-boundedness'', 
{\em J. Graph Theory} {\bf 95} (2020), 473--504.

\bibitem{poly2} A. Scott, P. Seymour and S. Spirkl, 
``Polynomial bounds for chromatic number. II. Excluding a star forest'', 
{\em J. Graph Theory}, {\bf 101} (2022), 318--322,
{\tt arXiv:2107.11780}.


\bibitem{spencer}
J. Spencer,
``Asymptotic lower bounds for Ramsey functions'', 
{\em Discrete math.} {\bf 20} (1977), 69--76.

\bibitem{spirkl}
S. Spirkl,
{\em Cliques, Stable Sets, and Coloring in Graphs with Forbidden Induced Subgraphs},
Ph.D. thesis, Princeton University (2018).


\bibitem{sumner}
D.P. Sumner, ``Subtrees of a graph and chromatic number'', in
{\em The Theory and Applications of Graphs}, (G. Chartrand, ed.),
John Wiley \& Sons, New York (1981), 557--576.

\bibitem{turan}
P. Tur\'{a}n, ``"On an extremal problem in graph theory'', {\em Matematikai \'{e}s Fizikai Lapok} (in Hungarian), {\bf 48} (1941), 436--452.
\end{thebibliography}
\end{document}